\documentclass{amsart}

\title[Pure eigenstates]{Pure eigenstates for the sum of generators of the free group}

\author[W. Paschke]{William L. Paschke}
\address{Department of Mathematics \\
University of Kansas \\
Lawrence, KS 66045-2142}
\email{paschke@math.ukans.edu}

\newtheorem{thm}{Theorem}[section] 
\newtheorem{prop}[thm]{Proposition} 
\newtheorem{lem}[thm]{Lemma} 
\newtheorem{rem}[thm]{Remark}

\begin{document}

\begin{abstract} 
We consider certain positive definite functions on a finitely
generated free group G  that are defined with respect to a given basis in terms of word length
and the number of negative-to-positive generator exponent switches. Some of these functions are
eigenfunctions for right convolution by the sum of the generators, and give rise to irreducible
unitary representations of G. We show that any state of the reduced C*-algebra of G whose left
kernel contains a polynomial in one of the generators must factor through the conditional
expectation on the C*-subalgebra generated by that generator. Our results lend some support to
the conjecture that an element of the complex group algebra of G can lie in the left kernel of
only finitely many pure states of the reduced C*-algebra of G.
\end{abstract}

\

\maketitle
  
\section{Introduction}\label{intro} 

Let $G$ be the free group on $n$ generators $u_1, u_2, \ldots u_n,$ where $n \geq 2.$ We
will regard the complex group algebra ${\mathbb C}G$ variously as a subalgebra of the reduced
C*-algebra $C^*_r(G)$ (the operator norm closure of the image of ${\mathbb C}G$ under the left
regular representation of $G$ on $\ell^2(G)$) and as a subalgebra of the full group
C*-algebra $C^*(G)$ (the completion of ${\mathbb C}G$ in the norm obtained by taking the supremum
over all unitary representations of $G$). Positive definite functions on $G$ will be
thought of as extended linearly to positive functionals on the *-algebra ${\mathbb C}G$ and thence
to positive linear functionals on $C^*(G)$ or, if appropriate, on $C^*_r(G);$ a positive
functional on ${\mathbb C}G$ that extends positively to $C^*_r(G)$ will be called reduced. The term
state refers to positive definite functions or positive linear functionals which take the
value 1 at the identity of $G$. Pure states are extreme points of the set of states of the
the relevant *-algebra; they give rise by a familiar construction to irreducible unitary
representations of $G$. For a scalar $\lambda$ and an algebra element $x$, a
$\lambda-$eigenstate of $x$ is a state with $x - \lambda$ in its left kernel, that is, a
state annihilating $(x-\lambda)^*(x - \lambda).$

The work reported on below is motivated in large part by a conjecture about how
${\mathbb C}G$ fits into $C^*_r(G)$, namely that each nonzero element of ${\mathbb C}G$ belongs to the left 
kernel of only finitely many pure states of $C^*_r(G)$. In other words, we surmise that
for each nonzero $y$ in ${\mathbb C}G$, the (convex) set of states on $G$ weakly associated to the
left regular representation and annihilating $y^*y$ has only finitely many extreme points.
We will show this when $y$ is the sum of the generators minus a scalar of modulus
$\sqrt{n}$, in which case the set of states in question is a singleton, and also when $y$
is a polynomial in one of the generators. Our results on $\lambda-$eigenstates of $u_1 +
u_2 + \ldots + u_n$ for $|\lambda| < \sqrt{n}$ point strongly in the direction of the
conjecture, but leave open the faint possibility that there might be an essential
difference between the interior and the boundary of the spectrum of $u_1 + u_2 + \ldots +
u_n$ as an operator on $\ell^2(G).$

We remark that such operators can have no eigenvectors in $\ell^2(G).$ This follows from a
result of P.A. Linnell \cite{Linnell} stating that for many torsion-free groups, including free
groups, the $\ell^2-$kernel of a matrix with entries in the group algebra must have integer von
Neumann dimension. In particular, no nonzero element of the group algebra can convolve a
nonzero $\ell^2-$function to zero.

Whatever the fate of the conjecture, we hope to persuade the reader that it is edifying to
look for pure eigenstates, especially but not exclusively those that come from
$C^*_r(G)$, of particular elements of ${\mathbb C}G$.  An example of what can come from this 
pursuit under favorable circumstances is the program recounted in \cite{FigaTalamancaPicardello}
of harmonic analysis on $G$ from the point of view of  radial functions, in which pure
eigenstates for the symmetrized  sum of generators
$$u_1 + u_1^{-1} + u_2 + u_2^{-1} + \ldots + u_n + u_n^{-1}$$ play a central role.
(Uniqueness of the eigenstates for this element among radial functions is  established
fairly easily in \cite{FigaTalamancaPicardello}; the question of a whether a
$C^*_r(G)-$eigenstate for the symmetrized sum must be radial is open.) See also
\cite{FigaTalamancaSteger} for harmonic analysis based on arbitrary linear combinations of the $u_j
+ u_j^{-1}$ 's.

Our treatment of eigenstates for the unsymmetrized sum in Section \ref{eigenstates} below begins
with the definition of a certain family of functions $\phi$ on $G$. We show by calculating  matrix
eigenvalues that each such $\phi$ is so to speak positive definite over the positive 
semigroup of $G$. This yields a Hilbert space on which the $u_j$'s act isometrically. An 
appropriate dilation then yields a unitary representation of $G$ on a larger Hilbert space
from which $\phi$ can be recovered by composing with a vector state. Using the result of
Linnell  mentioned above, we show that $u_1 + u_2 + \ldots + u_n$ has exactly one
eigenvalue (depending  on which $\phi$ one starts with) in this representation, with a
one-dimensional eigenspace. It follows from this that different $\phi$'s give rise to
unitarily inequivalent irreducible unitary  representations of $G$. In Section \ref{more} we
consider these functions in the context of a somewhat larger set of states, defined, like
the original $\phi$'s, in terms of word length and the number of negative-to-positive
generator exponent changes. We determine which of these are states of the reduced
C*-algebra. We also indicate how the irreducible representations in Section \ref{eigenstates}
(for an appropriate range of spectral values) can be realized in terms of the action of $G$ on
its combinatorial boundary; the measures on the boundary that we consider are rather like
those treated by G. Kuhn and T. Steger in \cite{KuhnSteger}. We prove uniqueness of the
$C^*_r(G)-$eigenstate of the sum of the generators for eigenvalues of modulus $\sqrt{n}$ in
Section \ref{paucity}, and show also that states of
$C^*_r(G)$ whose left kernel contains a given polynomial in a generator $u_j$ must factor
through the conditional expectation on the C*-algebra generated by $u_j$. When we identify
this C*-algebra with the algebra of continuous functions on the unit circle, the relevant
pure states are point evaluations at zeros of modulus 1 of the given polynomial preceded
by the conditional expectation; hence, there are only finitely many such pure states.

\

Heartfelt thanks are owed the referee of this paper, who patiently pointed out numerous {\em
faux pas}, and suggested that the rather cumbersome argument originally provided for Theorem
\ref{poly}  be replaced by the straightforward proof that now appears.

\

\section{The eigenstates and their representations}\label{eigenstates}

\

In seeking eigenstates, reduced or not, for $u_1 + u_2 + \ldots + u_n$, we may confine our
attention to spectral values in the interval $[0,n]$. This is because for each complex
$z$ of modulus 1 there are automorphisms of both $C^*_r(G)$ and $C^*(G)$ sending each
$u_j$ to $z u_j.$ It is convenient to divide by $n$ and use $[0,1]$ for the parameter
interval. Thus $a$ in the unit interval corresponds to the spectral value $na$. For such $a$,
define $\phi_a$ on $G$ by 
$$\phi_a(s) = a^{|s| - 2 \gamma(s)} \left( \frac{na^2 - 1}{n-1} \right)^{\gamma(s)} ,$$
where $|s|$ is the length of $s$ as a reduced word in the given generators and their
inverses, and $\gamma(s)$ is the number of negative-to-positive generator exponent
changes in $s$. (Thus for example $\gamma(u_1^{-2}u_2^3u_1^{-1}) = 1.$) We will often
write $\phi$ instead of $\phi_a$. In case $a$ is 0 or $1/\sqrt{n}$, we take $0^0$ to be 1.
Let us check that
$$\sum_{j=1}^n \phi(su_j) = na \phi(s)$$
for every $s$ in $G$. Since $\phi(1) = 1$ and $\phi(u_j) = a$ for each $j$, this is true
when $s = 1.$ If $s$ ends in $u_i^{-1}$ for some $i$, then $|su_j| = |s| + 1$ and
$\gamma(su_j) = \gamma(s) + 1$ for $j \not= i,$ while $|s u_i| = |s| - 1$ and
$\gamma(su_i) = \gamma(s),$ so
$$\sum_{j=1}^n \phi(su_j) = \left(\frac{n-1}{a} \cdot \frac{na^2 -1}{n-1} +
\frac{1}{a}\right)
\phi(s) = na \phi(s).$$
(For the case $a = 0,$ take limits in this calculation.) If $s$ ends in $u_i$ for
some $i$, then $\phi(su_j) = a \phi(s)$ for every $j$, and again the formula holds.
Notice, incidentally, that $|s^{-1}| = |s|$ and $\gamma(s^{-1}) = \gamma(s)$ for every
$s$ in $G$, so $\phi$ is selfadjoint.

Our main project in this section is to show that $\phi$ is positive definite and that the
unitary representation of $G$ to which it gives rise is irreducible. The case $a = 0$,
alas, requires somewhat special treatment, so we will assume until further notice that $0
< a \leq 1.$ To save  space, let us write
$$b = \frac{na^2 - 1}{n-1} \ .$$

Let $G^+$ be the unital semigroup  in $G$ generated by $u_1, u_2, \ldots, u_n$,
and let $G^+_k$ be the set of group elements in $G^+$ of length $k$. For $k = 1, 2,
\ldots,$ let $A_k$ be the $n^k \times n^k$ matrix with entries indexed by $G^+_k \times
G^+_k$ whose $(s,t)$-entry is $\phi(s^{-1}t).$ Since $\phi$ is real-valued and satisfies
$\phi(s) = \phi(s^{-1}),$ the matrices $A_k$ are hermitian.

\

\begin{lem}\label{pd} The matrix $A_k$ is positive semidefinite for $k = 1, 2, \ldots \ .$
\end{lem}

\

\noindent{\bf Proof:} Notice that $A_1$ is the $n \times n$ matrix with $1$'s on the diagonal and
$b$ in every off-diagonal position. The spectrum of $A_1$ is easily seen to be
$$\{1 - b, \ 1 + (n-1) b\} \ ,$$ so since $-(n-1)^{-1} < b \leq 1,$ we have 
$A_1 \geq 0.$ For the inductive step, regard $G^+_{k+1}$ as the disjoint union of $n$ 
copies of $G^+_k$ by writing
$$G^+_{k+1} = u_1 G^+_k \cup u_2 G^+_k \cup \dots \cup u_n G^+_k \ .$$ We can then write
$A_{k+1}$ in terms of $A_k$ as an $n \times n$ matrix of $n^k \times n^k$ matrices.
Namely, $A_{k+1}$ has $A_k$ in each of the $n$ diagonal blocks, and all of its
off-diagonal blocks are $a^{2 k} b$ times the $n^k \times n^k$ matrix $E_k$ with all
entries equal to $1$. This is because $\phi(s^{-1}u_i^{-1}u_j t),$ for $s, t$ in $G^+_k$,
is $\phi(s^{-1}t)$ if $i = j$ and $a^{2k} b$ if $i \not= j$. It follows by induction on
$k$ that the entries in every row of $A_k$ sum to the common value
$$ 1 + (n-1) b + (n-1) n a^2 b + \ldots + (n-1) n^{k-1} a^{2k - 2} b = (n a^2)^k \ .$$ Let
$\lambda$ be an eigenvalue of $A_{k+1}$ which is not an eigenvalue of $A_k$.  The
corresponding eigenvector is an $n$-tuple $(\xi_1, \xi_2, \ldots, \xi_n)$ of vectors 
$\xi_j$ with entries indexed by $G^+_k$ satisfying
$$(A_k - \lambda) \xi_i + a^{2k} b E_k (\sum_{j \not= i} \xi_j) = 0$$ for $i = 1, 2,
\ldots, n.$ The range of $E_k$, vectors with all entries the same, is invariant under
$A_k$ (since the latter has constant row sums). Because $A_k - \lambda$ is invertible, it
follows that each $\xi_i$ belongs to the range of $E_k$. Let $c_i$ denote the common value of the
entries of $\xi_i$. Then
$$((n a^2)^k - \lambda) c_i + n^k a^{2k} b \sum_{j \not= i} c_j = 0$$ for $i = 1, 2,
\ldots n,$ which is to say that $\lambda$ is an eigenvalue of $(n a^2)^k A_1,$ hence
nonnegative. If we know that $A_k \geq 0,$ then $A_{k+1} \geq 0.$

\

The Hilbert space $H$ of the representation we seek is constructed as follows. By
Lemma \ref{pd}, there is for each positive integer $k$  a finite dimensional complex inner
product space  $E_k$ spanned by vectors $\{\Delta_s : s \in G^+_k \}$ with inner
product $< \cdot, \cdot >$ satisfying $<\Delta_t , \Delta_s> = \phi(s^{-1}t).$ (We write 
$E_0$ for the one-dimensional inner product space spanned by the unit vector
$\Delta_1.$) Because  
$$(na)^{-1} \phi(\sum x u_i) = \phi(x) = (na)^{-1} \phi(\sum  u_i^{-1} x)$$
for all $x$ in ${\mathbb C}G$, we have an isometry from $E_k$ into $E_{k+1}$
for each $k$ sending $\Delta_s$  to $(na)^{-1} \sum_i \Delta_{su_i}.$ Let $H_0$ be the Hilbert space
inductive limit of the resulting tower $E_0 \rightarrow E_1  \rightarrow E_2
\rightarrow \ldots .$ Thus $H_0$ is the closed linear span of $\{\Delta_s : s \in
G^+ \}$, and these vectors satisfy 
$$\sum_i \Delta_{su_i} = na \Delta_s \ \ \mbox{and} \ \  <\Delta_t , \Delta_s> \ = \
\phi(s^{-1}t).$$  
Left multiplication by each generator $u_i$ gives rise to an isometry $V_i$ of 
$H_0$ into itself. Let $H_i'$ be the kernel of $V_i,$ in other words, the
orthogonal complement in $H_0$ of the range of $V_i.$ For each $i,$ let
$S^-_i$ be the subset of $G$ consisting of the reduced words ending in a negative power of
$u_i$, with natural orthonormal basis $\{\delta_s : s \in S^-_i \}$. The Hilbert space
$H$ is by definition
$$H = H_0 \oplus \bigoplus_{i = 1}^n \left(\ell^2 (S^-_i) \otimes H_i'
\right) .$$  
For each $i$, let $U_i$ be the unitary operator on $H$ that maps $H_0$ to $V_i H_0$ by $V_i$, maps
$\delta_{u_i^{-1}} \otimes H_i'$ to $H_i' =  H_0 \ominus V_iH_0$ by erasing the tensor, and maps
$\delta_s
\otimes \eta$ to $\delta_{u_i s} \otimes \eta$ for all other $s$ ending in a negative
generator power, and for $\eta$ in the appropriate space $H_j'$. Denote by $\pi$ the
unitary representation of $G$ on $H$ that takes $u_i$ to $U_i$.

We now show that $\phi(s) \ = \ <\pi(s) \Delta_1, \Delta_1>$ for all $s$ in $G$. This is
mostly a matter of decomposing $\pi(s) \Delta_1$ into orthogonal pieces
as in the definition of $H.$

\

\begin{lem}\label{orthog} For $i = 1, 2, \ldots, n$, and $s \in G^+ \setminus \{1\}$ not
beginning with $u_i$, the vectors
$$\Delta_1 - a \Delta_{u_i} \ \ \mbox{and} \ \ \Delta_s - a^{|s| - 1} b \Delta_{u_i}$$
both belong to $H_i'$, and hence
$$U_i^* \Delta_1 = a \Delta_1 + \delta_{u_i^{-1}} \otimes (\Delta_1 - a \Delta_{u_i}),$$
$$U_i^* \Delta_s = a^{|s| - 1} b \Delta_1 + \delta_{u_i^{-1}} \otimes  (\Delta_s - a^{|s|
- 1} b \Delta_{u_i}).$$
\end{lem}

\

\noindent{\bf Proof:} Take $t$ in $G^+.$ Calculating in $H_0,$ we have
$$< \Delta_1 - a \Delta_{u_i}, \, \Delta_{u_i t} > \ = \ \phi(t^{-1} u_i^{-1}) - a
\phi(t^{-1}) = a^{|t| + 1} - a a^{|t|} = 0,$$ and
$$<\Delta_s - a^{|s| - 1} b \Delta_{u_i}, \, \Delta_{u_i t} > \ = \ \phi(t^{-1} u_i^{-1} s)
- a^{|s| - 1} b \phi(t^{-1}) $$
$$= a^{|t| + 1 + |s| - 2} b - a^{|s|-1} b a^{|t|} = 0.$$  Thus $\Delta_1$ and $\Delta_s$
are written with respect to the orthogonal decomposition
$H_0 = U_i H_0 \oplus H_i'$ as
$$\Delta_1 = a \Delta_{u_i} + (\Delta_1 - a \Delta_{u_i}) \ , \  
\Delta_s =  a^{|s| - 1} b \Delta_{u_i} + (\Delta_s - a^{|s| - 1} b \Delta_{u_i}) \ ,$$ 
and the rest of the lemma follows by noticing that $U_i^* \Delta_{u_i} = \Delta_1$, while
$U_i^*\eta = \delta_{u_i^{-1}} \otimes \eta$ for $\eta$ in $H_i'.$

\

\begin{prop}\label{phi} $\phi(s) \ = \ <\pi(s) \Delta_1, \Delta_1>$ for all 
$s$ in $G$.
\end{prop}

\

\noindent{\bf Proof:} Write $\psi(s) \ = \ <\pi(s) \Delta_1, \Delta_1>$. The argument  that $\psi(s)
=\phi(s)$ is by induction on $|s|$. The cases $|s| = 0, |s| = 1$ are clear. For (part of)
the induction, consider $\psi(t u_i^{-1})$ versus $\psi(t)$, where
$t$ is a reduced word in $G$ that doesn't end in a positive power of $u_i$. We  have 
$$\pi(t) U_i^* \Delta_1 = \pi(t) (a \Delta_1 + \delta_{u_i^{-1}} \otimes \eta)$$ for
appropriate $\eta$ in $H_i'$ by Lemma \ref{orthog}. Furthermore, the definition of
$\pi$ and our assumption on $t$ ensure that 
$$\pi(t) (\delta_{u_i^{-1}} \otimes H_i') \subseteq  
\ell^2(S_i^-) \otimes H_i' \ .$$ Since $\Delta_1$ is orthogonal to the latter
subspace, we have $\psi(t u_i^{-1}) = a \psi(t).$ In the same situation, we can compare
$\psi(t u_i^{-1} u_j)$ with 
$\psi(t u_i^{-1})$ when $i \not= j$. Indeed, Lemma \ref{orthog} shows that for appropriate
$\xi$ in $H_i'$ we have
$$\pi(t) U_i^*U_j \Delta_1 = \pi(t) U_i^*\Delta_{u_j} = 
\pi(t) (b \Delta_1 + \delta_{u_i^{-1}} \otimes \xi) \ ,$$ and hence $\psi(t u_i^{-1} u_j) = b
\psi(t) = (b/a) \psi(t u_i^{-1})$. We must show as well that $\psi(s u_j) = a \psi(s)$
for every generator $u_j$ if $s$ ends in a positive generator power. It follows from Lemma
\ref{orthog} and the definition of $\pi$ that for such an $s$, we have $\pi(s^{-1}) \Delta_1 = c
\Delta_1 + \rho,$ where $\rho$ is orthogonal to $H_0$, and $c$ is a real constant
namely $c = \psi(s^{-1}) = \psi(s).$ (Start with $\Delta_1$ and apply each factor $u_i^{\pm}$
of $s^{-1}$ in succession. The last factor applied has exponent $-1$.) Thus
$$U_j^*\pi(s^{-1}) \Delta_1 = ca \Delta_1 +  \delta_{u_j^{-1}} \otimes \eta +  U_j^* \rho$$
for appropriate $\eta$ in $H_j'$. Since $H \ominus H_0$ is  invariant
under $U_j^*$, this makes
$\psi(u_j^{-1} s^{-1}) = ca = a \psi(s),$ and finally $\psi(s u_j) = a \psi(s)$. Thus,
when the length of a word is increased by 1 by non-cancelling right multiplication  by a
generator or its inverse, the value of $\psi$ is multiplied by $b/a$ or $a$ depending on
whether or not $\gamma$ increases by 1. This is the same rule that $\phi$ obeys, so
$\phi$ and $\psi$ must coincide.

\

The next proposition implies that $\pi$ is irreducible, and that different 
$a$'s in the interval $(0,1]$ give rise to unitarily inequivalent representations. 

\

\begin{prop}\label{onlyeigen} The only eigenvalue that $u_1 + \ldots + u_n$ has in the
representation $\pi$ is $na.$ The eigenspace consists of scalar multiples of $\Delta_1.$
\end{prop}

\

\noindent{\bf Proof:} We have already observed that $\sum_j \Delta_{su_j} = na \Delta_s$ for all
$s$ in $G^+$. In particular, $(\sum_j U_j - na) \Delta_1 = 0.$ 

Suppose that $\lambda$ is an eigenvalue for
$\sum_j U_j$ with eigenvector $\xi.$ We first claim that $\xi$ must belong to 
$H_0.$ Fix a generator $u_i$, and take $\eta$ in $H_i'$. Define
$f_0$ on $S^-_i$ (the set of reduced words in $G$ ending in a negative power of
$u_i$) by $f_0(s) = <\xi, \delta_s \otimes \eta>.$ Then $f_0 \in \ell^2(S^-_i)$ and
satisfies $\sum_j f_0(u_j^{-1} s) = \lambda f_0(s)$ for all $s$ in $S^-_i$. This is because
$$\lambda f_0(s) \ = \ <(\sum_j U_j) \xi, \delta_s \otimes \eta> \ = \ <\xi, \sum_j
(\delta_{u_j^{-1} s} \otimes \eta)> .$$ Pick a different generator $u_h$ and let $\sigma$ be the
automorphism of 
$G$ that interchanges $u_i$ and $u_h$, and fixes the other generators. Define
$f$ in $\ell^2(G)$ by setting $f(s) = f_0(s)$ for $s$ in $S^-_i$, and $f(s) =
-f_0(\sigma(s))$ for $s$ in $S^-_h$, and $f(s) = 0$ for all other $s$ in $G$. This makes
$\sum_j f(u_j^{-1} s) = \lambda f(s)$ for all $s$ in $G$. (Equality holds for
$s$ in $S^-_i \cup S^-_h$ by construction. In case $s = 1,$ the right-hand side is 0 and
the left-hand side is $f(u_i^{-1}) + f(u_h^{-1}),$ which is 0. For all other
$s$, we have $f(s) = 0$, and each summand on the left is 0.) In other words, when we let
$G$ act on $\ell^2(G)$ by the left regular representation, the $\ell^2$ function
$f$ belongs to the kernel of $u_1 + \ldots + u_n - \lambda$ (as an operator on $\ell^2(G)$). By
3.6 in \cite{Linnell} (the result of  Linnell mentioned in the introduction), this forces $f,$
and hence $f_0$, to vanish everywhere. Since
$\eta$ in
$H_i'$ was arbitrary, we have  shown that $\xi$ is orthogonal to each summand
$\ell^2(S^-_i) \otimes H_i',$ which means $\xi \in H_0.$

Consider now $<\xi, \Delta_s>$ for $s$ in $G^+$. For such $s$ and any generator $u_j$, we
have
$$\lambda <\xi, \Delta_{u_j s}> \ = \ <\xi, (U_1^* + \ldots + U_n^*) \Delta_{u_j s}> .$$
Because $\xi \in H_0,$ we see by applying Lemma \ref{orthog} to the right-hand side that
$$(*) \ \ \  \lambda <\xi, \Delta_{u_j s}> \ = \ <\xi, \Delta_s> +  (n-1) a^{|s|}b <\xi,
\Delta_1> \ .$$ With $s = 1,$ this becomes
$$(**) \ \ \ \lambda <\xi, \Delta_{u_j}> \ = \ (1 + (n-1) b) <\xi, \Delta_1> \ = \ na^2
<\xi, \Delta_1> \ .$$ Since $\Delta_{u_1} + \ldots + \Delta_{u_n} = na \Delta_1,$ when we
sum on $j$ in (**) we obtain
$$ \lambda na <\xi, \Delta_1> \ = \ n^2 a^2 <\xi, \Delta_1>.$$ In case $<\xi,
\Delta_1> \not= 0,$ we have $\lambda = na,$ and when we replace  $\xi$ by
$\xi' = \xi - <\xi, \Delta_1> \Delta_1$ (a vector in the $na$-eigenspace meet 
$H_0$ orthogonal to $\Delta_1$), equation (*) above becomes
$$ na <\xi', \Delta_{u_j s}> \ = \ <\xi', \Delta_s>$$ for all generators $u_j$ and all $s$
in $G^+$. The right-hand side is zero when $s = 1,$ so by induction $\xi'$ is orthogonal to
all the $\Delta_s$'s, so $\xi' = 0,$ finishing the proof in this case. If, on the other
hand, $<\xi, \Delta_1> = 0,$ the same argument using (*) shows $\xi = 0$ if $\lambda \not=
0,$ while if $\lambda$ and $<\xi, \Delta_1>$ are both 0, then (*) simply says that $<\xi,
\Delta_s> = 0$ for all $s$ in $G^+$, and again $\xi = 0.$

\

We turn now to the hitherto excluded case $a = 0,$ which is best argued separately.

\

\begin{prop}\label{zerocase} Let $b = -1/(n-1),$ and define $\phi$ on $G$ by 
$\phi(s) = b^{\gamma(s)}$ if $|s| = 2 \gamma(s)$, and $\phi(s) = 0$ otherwise. Then
$\phi$ is positive definite. Furthermore, zero is the only eigenvalue of $u_1 + \ldots +
u_n$ in the cyclic unitary representation of $G$ on Hilbert space to which $\phi$ gives
rise, and the eigenspace is one-dimensional.
\end{prop}

\

\noindent{\bf Proof:} The reason $\phi$ is positive definite is that it is the pointwise limit of
functions that are positive definite by Proposition \ref{phi}. Looking
at the part of the corresponding representation space spanned by vectors from
$G^+$, we see that there is a Hilbert space
$H_0$ which is the closed linear span of unit vectors $\Delta_s$ ($s$ in $G^+$)
satisfying 
$$\sum_i \Delta_{su_i} = 0 \ \ \mbox{and} \ \  <\Delta_t , \Delta_s> \ = \
\phi(s^{-1}t).$$ (Notice that this makes $\Delta_s$ and $\Delta_t$ orthogonal if either
$s$ or $t$ has length greater than 1.)   Using the same notation as in the construction
immediately before Lemma \ref{orthog}, we obtain the  Hilbert space 
$$H = H_0 \oplus \bigoplus_{i = 1}^n \left(\ell^2 (S^-_i) \otimes H_i'
\right) $$ and unitary  operators $U_1, \ldots, U_n$ on $H$ ---  with $U_i$ taking
$\delta_{u_i^{-1}} \otimes \eta$ to $\eta$ for $\eta$ in 
$H_i' = H_0 \ominus U_i H_0$, and so forth. This apparatus gives rise 
in turn to a unitary representation  $\pi$ of $G$ on $H$ with cyclic vector 
$\Delta_1$. 

Define $\psi$ on $G$ by $\psi(s) = <\pi(s) \Delta_1, \Delta_1>.$ By construction,
$\psi$ and $\phi$ coincide on $G^+$; they are both 1 at 1 and vanish on the rest of
$G^+$. A reduced word in $G \setminus G^+$ can be written in the form $s u_j^{-1} v,$
where $v$ (possibly empty) belongs to $G^+$ and does not begin with a positive power of $u_j$,
and $s$ is a  reduced word (possibly empty) not ending in a positive power of $u_j$. The
definition of $\phi$ entails that  $\Delta_v \in H_j'$ if $v = 1$ or if $|v| > 1.$ In case
$v = u_i$ for some $i \not=  j$, the orthogonal projection of $\Delta_v$ on $H_j'$ is 
$\Delta_{u_i} - b \Delta_{u_j}$. Thus
$$\pi(u_j^{-1} v) \Delta_1 = \left\{ \begin{array}{cl}  b \Delta_1 + \delta_{u_j^{-1}}
\otimes (\Delta_{u_i} - b \Delta_{u_j})  & 
 \ \ \ v = u_i \ \mbox{for} \ \ i \not= j \\
\delta_{u_j^{-1}} \otimes \Delta_v &  \ \ \ v \in G^+ \setminus \{u_1, \ldots, u_n\}
\end{array} \right. \ \ .$$ This makes
$$\psi(s u_j^{-1} v) = \left\{ \begin{array}{cl}  b \psi(s)  & 
\ \ \ v = u_i \ \mbox{for} \ \ i \not= j \\ 0 & \ \ \ v \in G^+ \setminus \{u_1, \ldots,
u_n\}
\end{array} \right. \ \ .$$ The same holds for $\phi$, whence it follows that $\phi$ and
$\psi$ coincide on all of $G$. 

Suppose now that $\lambda$ is an eigenvalue of $U_1 + \ldots + U_n$ with eigenvector
$\xi.$ We must show that $\lambda = 0$ and $\xi$ is a scalar multiple of $\Delta_1.$ That
$\xi$ belongs to $H_0$ follows exactly as in the proof of Proposition \ref{onlyeigen}. For
any $s$ in $G^+$ and any generator $u_j$ we have
$$\lambda <\xi, \Delta_{u_j s}> \ = \ <\xi, \sum_i U_i^* \Delta_{u_j s} > \ = \ \left\{
\begin{array}{cl} <\xi, \Delta_s> & \ \ s \in G^+ \setminus \{1\} \\ 0 & \ \ s = 1
\end{array} \right. \ \ .$$ This is because in the first case, $U_i^* \Delta_{u_j s}$ is orthogonal
to $H_0$, hence to $\xi$, for $i \not= j$, and in the second case, the right-hand side is 
$$<\xi, \Delta_1 + b(n-1) \Delta_1> \ = \ <\xi, \Delta_1 - \Delta_1> \ .$$ 
If $\lambda$ is different from zero, we get $<\xi, \Delta_{u_j}> = 0$ for all $j$ and then 
$<\xi, \Delta_s> = 0$ for all other $s$ in $G^+ \setminus \{1\},$ forcing $\xi$ to be a
multiple of $\Delta_1,$ so $\xi = 0.$ On the other hand, $\lambda = 0$ forces
$<\xi, \Delta_s> = 0$ for $s$ in $G^+ \setminus \{1\},$ so $\xi$ in this case must be a
multiple of $\Delta_1.$

\

We now summarize the results of this section.

\

\begin{thm}\label{summarize} Let $G$ be the free group on the generators $u_1, \ldots,
u_n$, where $n \geq 2.$ Let $| \cdot |$ denote the corresponding length function on $G$, 
and let $\gamma$ be the function on $G$ that counts the number of negative-to-positive
exponent changes, from left to right, in a reduced word in these generators. Given $0 \leq a
\leq 1,$ define $\phi$ on $G$ by
$$\phi(s) = a^{|s| - 2 \gamma(s)}\left(\frac{n a^2 - 1}{n - 1}\right)^{\gamma(s)}.$$ Then
$\phi$ is a pure na-eigenstate of $G$, the unitary representation of $G$ on Hilbert space
to which $\phi$ gives rise is irreducible, and different values of  $a$ yield unitarily
inequivalent representations.
\end{thm}

\

\noindent{\bf Proof:} We have exhibited a representation $\pi$ of $G$ on a Hilbert space $H$
with cyclic vector $\Delta_1$ such that $\phi(s) = <\pi(s) \Delta_1, \Delta_1>$ for
$s$ in $G$. We have further shown that the kernel of $\pi(u_1) + \ldots + \pi(u_n) - na$
consists of scalar multiples of $\Delta_1.$ Thus, projection on the subspace spanned by
$\Delta_1$ belongs to the double commutant of $\pi(G).$ By cyclicity of $\Delta_1,$ it
follows that the double commutant must include all bounded operators on $H$, which
means that $\pi$ is irreducible. Finally, the point spectrum of 
$\pi(u_1) + \ldots + \pi(u_n)$, namely $\{na\}$ by Propositions \ref{onlyeigen} and
\ref{zerocase}, is invariant under unitary transformation of $\pi.$

\

For other spectral values, rotate Theorem \ref{summarize}.  More precisely, to multiply
the spectral value  $na$ by a phase $e^{i \theta}$ and thereby cover the entire spectrum of $u_1
+ \ldots + u_n$ relative to $C^*(G)$, simply compose $\phi$ with the automorphism of ${\mathbb C}G$
that multiplies each $u_j$ by $e^{i \theta}$. Let  $\tau : G \rightarrow Z$ be the homomorphism
sending each $u_j$ to 1. Then
$$ s \mapsto e^{i \tau(s) \theta} a^{|s| - 2 \gamma(s)}\left(\frac{n a^2 - 1}{n -
1}\right)^{\gamma(s)}$$
is a pure $na e^{i \theta}-$ eigenstate of $G$, and different spectral values give rise to
unitarily inequivalent irreducible representations. Multiplying by $-1$, in particular, we
see that the interval [0,1] in Theorem \ref{summarize} can be replaced by the interval [-1,1].

\

\section{More states; reduced states}\label{more}

\

For $a$ in the interval $[-1,1]$, the function $s \mapsto a^{|s| - 2 \gamma(s)}$ is
positive definite on $G$. This follows from \cite{DeMicheleFigaTalamanca} after a simple change
of generators. Namely, consider the ``$u_1-$length'' function $| \cdot |_1$ on $G$ that adds
the absolute values of all exponents of $u_1$ in a reduced word. Notice that $|st|_1 =
|s|_1 + |t|_1$ whenever there is no cancellation in multiplying $s$ and $t$, so $s
\mapsto a^{|s|_1}$ is positive definite by Theorem 1 of \cite{DeMicheleFigaTalamanca}. To
apply this to our situation, consider the automorphism $\beta$ of $G$ that fixes
$u_1$ and takes $u_j$ to $u_1 u_j$ for $j = 2, \ldots, n;$ one checks easily that
$|\beta(s)|_1 = |s| - 2 \gamma(s)$ for $s$ in $G$. 

\

Positive-definiteness of $s \mapsto a^{|s|_1}$ is all we need for the sequel, but it will
not take us too far afield to look at the corresponding unitary representations in the
spirit of the partial tensor product construction used in the  previous section.  Fix $a$
in $[-1,1]$ and consider the positive definite function $k \mapsto a^{|k|}$ on the group
of integers. (The reason this is positive definite is that
$$1 + \sum_{k=1}^\infty a^k (e^{i k \theta} + e^{-i k \theta}) = 
\frac{1- a^2}{1 + a^2 - 2 a \cos \theta}$$ for $-1 < a < 1.$) Let $\pi_0$ be the
associated unitary representation on the Hilbert space $H_0$, with cyclic vector
$\Delta_1,$ thought of as a representation of the subgroup of $G$ generated by $u_1.$
Thus, $<\pi_0(u_1)^k \Delta_1, \Delta_1> = a^{|k|}$ for all integers $k$. Let $H_0'$
be the subspace of $H_0$ orthogonal to
$\Delta_1;$ it is more or less immediate that the vectors 
$\pi_0(u_1^k) \Delta_1 - a^{|k|} \Delta_1$ for $k \not= 0$ span a dense subspace of 
$H_0'$. Let $S$ be the subset of $G$ consisting of all nonempty reduced words not
ending in a power of $u_1,$ that is, all reduced words ending in a power of $u_j$ for some
$j > 1.$ Define the Hilbert space $H$ by
$$H = H_0 \oplus \left(\ell^2(S) \otimes H_0' \right) \ .$$ 
Let $U_1$ be the unitary operator on $H$ such that $U_1 \xi = \pi_0(u_1) \xi$ for $\xi$ in
$H_0$ and $U_1 (\delta_s \otimes \eta) = \delta_{u_1 s} \otimes \eta$ for $s$ in $S$
and $\eta$ in $H_0'$. Define unitaries $U_2, \ldots, U_n$ by:
$$U_j \Delta_1 = \Delta_1 \ ;  \ \ U_j \eta = \delta_{u_j} \otimes \eta \ ; $$
$$U_j (\delta_{u_j^{-1}} \otimes \eta) = \eta \ ; \ \ \ \mbox{and} \ \ \ U_j (\delta_s
\otimes \eta) = \delta_{u_j s} \otimes \eta$$ for $\eta$ in $H_0'$ and $s$ in $S
\setminus \{u_j^{-1}\}$. Let $\pi$ be the representation of $G$ on $H$ that takes
$u_i$ to $U_i$ for $i = 1, \ldots, n.$

\

\begin{prop}\label{pi} The representation $\pi$ satisfies 
$<\pi(s) \Delta_1, \Delta_1> = a^{|s|_1}$ for all $s$ in $G$, and is irreducible. Different
values of $a$ in $[-1,1]$ give rise to unitarily inequivalent representations.
\end{prop}

\

\noindent{\bf Proof:} Notice that $\pi(s) H_0' = \delta_s \otimes H_0'$ for $s$ in
$S$, so in particular $\pi(s) H_0'$ is orthogonal to  $\Delta_1$. It follows that
$$<\pi(s u_1^k) \Delta_1, \Delta_1> \ = \ <\pi(s) (\pi(u_1^k) \Delta_1 - a^{|k|} \Delta_1 +
a^{|k|} \Delta_1) , \Delta_1>$$
$$ =  a^{|k|} <\pi(s) \Delta_1, \Delta_1>$$ for $s$ in $S$ and any integer $k$, because
$\pi(u_1^k) \Delta_1 - a^{|k|} \Delta_1 \in  H_0'$. It is furthermore plain that
multiplying $s$ in $G$ on the left  by a power of $u_j$ for $j \geq 2$ does not change the value
of $<\pi(s) \Delta_1, \Delta_1>.$ Our formula for $<\pi(\cdot) \Delta_1, \Delta_1>$ follows. 

To distinguish between different values of $a$, fix sequences $\{s_k\}, \{t_k\}$ in the
subgroup of $G$ generated by $u_2, \ldots, u_n$ such that $|s_k| \rightarrow \infty$ and 
$|t_k| \rightarrow \infty$ as $k \rightarrow \infty.$ We claim that the sequence $\{\pi(s_k
u_1 t_k)\}$ converges in the weak operator topology to $a$ times the orthogonal projection
on $\Delta_1.$ To see this, let $\eta = U_1 \Delta_1 - a \Delta_1,$ so $\eta \in H_0',$
and $\pi(s_k) \eta = \delta_{s_k} \otimes \eta.$ This means that 
$\pi(s_k) \eta \rightarrow 0$ weakly as $k \rightarrow \infty$ and hence
$$\pi(s_k u_1 t_k) \Delta_1 = \pi(s_k u_1) \Delta_1 = \pi(s_k) (a \Delta_1 + \eta) 
\rightarrow a \Delta_1$$ weakly. For any $\eta$ in $H_0'$, we have $\pi(s_k u_1
t_k) \eta =
\delta_{s_k u_1 t_k} \otimes \eta$ for all $k$, and for any $s$ in $S$, we have
$\pi(s_k u_1 t_k) \delta_s \otimes \eta = \delta_{s_k u_1 t_k s} \otimes \eta$ for
sufficiently large $k$. It follows that $\pi(s_k u_1 t_k) \xi \rightarrow 0$  weakly for
any $\xi$ in $H$ orthogonal to $\Delta_1,$ which establishes our claim. The property
that all sequences of this type have weak-operator limit a certain scalar times a one-dimensional
projection is of  course invariant under unitary transformation of $\pi.$

If $a \not = 0,$ the argument just given also shows that $\pi$ is irreducible, since $\Delta_1$
is a cyclic vector. In case $a = 0,$ the representation $\pi$ is the one that comes from the
left action of $G$ on its quotient by the subgroup generated by $u_2, \ldots, u_n,$ and
irreducibility can be proved directly in several ways.

\

We now consider a two-parameter family of states on $G$ defined in terms of $| \cdot |$
and $\gamma.$ For real $a$ and $b$, define $\psi_{a,b}$ on $G$ by 
$$\psi_{a,b}(s) = a^{|s| - 2 \gamma(s)} b^{\gamma(s)} \ .$$

\

\begin{prop}\label{psi} The function $\psi_{a,b}$ is positive definite if and only if
$$-1 \leq a \leq 1 \ \ \mbox{and} \ \ \frac{na^2 - 1}{n-1} \leq b \leq 1 .$$
\end{prop}

\

\noindent{\bf Proof:} Assume that $\psi_{a,b}$ is positive definite. The inequalities
$|a| \leq 1$ and $b \leq 1$ are immediate because $\psi_{a,b}$ takes the values $1, a,$
and $b$ at $1, u_1,$ and $u_1^{-1}u_2$ respectively. The quadratic inequality comes
from the observation that
$$\psi_{a,b} \left( (na - \sum_i u_i^{-1}) (na - \sum_i u_i) \right) = $$
$$n^2 a^2 - 2 n a \sum_i \psi_{a,b}(u_i) + \sum_{i,j} \psi_{a,b}(u_i^{-1}u_j)
\ = \ -n^2 a^2 + n + n (n-1) b.$$

Suppose conversely that the indicated inequalities hold. It is then straightforward to find
$\alpha$ and $r$ with 
$$-1 \leq r \leq 1, \ 0 \leq \alpha \leq 1,  \ a = r \alpha, \ \mbox{and} \ 
b = \frac{n \alpha^2 - 1}{n - 1} \ .$$
It follows from Proposition \ref{pi} that $\psi_{r,1}$ is positive definite, while $\psi_{\alpha,
b}$ is positive definite by Proposition \ref{phi}, since $\psi_{\alpha, b} = \phi_\alpha.$ We
have
$\psi_{a,b} = \psi_{r,1} \psi_{\alpha, b},$ so $\psi_{a,b}$ is positive definite. 

\

We next sort out which values of $(a,b)$ in the region specified by Proposition \ref{psi} give
rise to reduced states of $G$, that is, to states whose associated unitary representation is
weakly contained in the regular representation.

\

\begin{prop}\label{reduced}  For $a, b$ as in Proposition \ref{psi}, the state
$\psi_{a,b}$ is reduced if and only if
$$ b \leq \frac{1 - na^2}{n - 1} .$$
\end{prop}

\

\noindent{\bf Proof:} Fix $a$ and $b$ satisfying the inequalities in Proposition \ref{psi}, and
write $\psi = \psi_{a,b}.$ By Theorem 3.1 in \cite{Haagerup}, whether or not $\psi$ is reduced
depends on whether the radius of convergence of the power series with order $k$ coefficient
$$C_k \equiv \sum_{|s| = k} |\psi(s)|^2$$  has radius of convergence at least 1, or
less than 1. We can obtain the $C_k$'s explicitly in the present instance by solving a pair
of one-step linear difference equations. For $k \geq 1,$ let
$$A_k = \sum \{|\psi(s)|^2 : |s| = k, s \ \mbox{ends in a positive generator power} \ \}$$ and
let $B_k$ be the corresponding sum over words ending in a negative generator power, so $C_k =
A_k + B_k .$ Notice that $A_1 = B_1 = na^2.$

If $a \not=0$, we have 
$$A_{k+1} = n a^2 A_k + (n-1) \frac{b^2}{a^2} B_k$$
$$B_{k+1} = (n-1)a^2 A_k + n a^2 B_k \ .$$ (For the first equation, take $s$ ending in a
positive generator power with $|s| = k+1$. Then either $s = tu_iu_j$ for some $j$, with $n$
choices for $i$ and
$\psi(s) = a \psi(tu_i),$ or $s = tu_i^{-1}u_j,$ with $n-1$ choices for $i$ and
$\psi(s) = (b/a) \psi(tu_i^{-1})$. The second equation is proved similarly.) If furthermore
$b \not= 0,$ the recurrence has two distinct eigenvalues, namely 
$$\lambda_+ = na^2 + (n-1)b \ \ \ \mbox{and} \ \ \ \lambda_- = na^2 - (n-1)b$$ with
eigenvectors $(\pm b, a^2)$. The solution may be written
$$A_k = \frac{n}{2} (a^2 + b) \lambda_+^{k-1} + \frac{n}{2} (a^2 - b) \lambda_-^{k-1}$$
$$B_k = \frac{n}{2}(a^2 + a^4/b) \lambda_+^{k-1} + \frac{n}{2}(a^2 - a^4/b) \lambda_-^{k-1}$$
By continuity, this holds as well when $a = 0,$ so the solution above covers the entire case
$b \not = 0.$ When $b = 0$, the coefficient matrix is lower triangular and one  calculates
directly that
$$A_k = \lambda^k \ \ \mbox{and} \ \ B_k = \frac{kn - k + 1}{n}  \lambda^k,$$ where $\lambda$ is
the common value of $\lambda_+$ and $\lambda_-$, namely $na^2.$

If 
$$b > \frac{1-na^2}{n-1} \ ,$$
we have $\lambda_+ > 1,$ and hence the power series with coefficients
$C_k$ has radius of convergence less than 1. If 
$$b \leq \frac{1-na^2}{n-1} \ ,$$ then also
$|a| \leq \sqrt{n}$ and  
$$|b| \leq \frac{1-na^2}{n-1} \ ,$$ which makes $|\lambda_+| \leq 1$ and
$|\lambda_-| \leq 1.$ In this case, then, the power series has radius of convergence at least 1.
  
\

\begin{rem}\label{contained} (a) If $a$ and $b$ satisfy $|b| < (1 - na^2)/(n-1),$ then $\sum C_k
< \infty,$ that is $\psi_{a,b} \in \ell^2(G),$ so the associated representation is contained
(not just weakly contained) in the regular representation. In particular, $\psi_{a,b}$ cannot
be a pure state in this situation.

(b) The proposition above shows that the state $\phi_a$ considered in Section \ref{eigenstates}
is reduced if and only if $|a| \leq 1/\sqrt{n}.$
\end{rem}

\

It is opportune to remark here as well that the spectral radius of $u_1 + \ldots + u_n$ as an
operator on $\ell^2(G)$ is indeed $\sqrt{n},$ as \ref{contained}(b) suggests. (The spectral radius is
at least $\sqrt{n}$ by \ref{contained}(b), and in the other direction, the operator norm of the
$k$th power is at most $(k+1)n^{k/2}$ by Haagerup's inequality \cite{Haagerup}.)

\

At least for spectral values of modulus strictly between 0 and $\sqrt{n}$,  the reduced
eigenstates in Section \ref{eigenstates} for the sum of the generators can be  obtained via the
action of $G$ on its combinatorial boundary. The compressed  account that follows is in the
spirit of G. Kuhn and T. Steger \cite{KuhnSteger}; the boundary representations considered here
differ only slightly from  the ones they treat. As in \cite{KuhnSteger}, let $\Omega$ denote the
set of all (one-way)  infinite reduced words in the generators $u_j$ and their inverses. When
given the product topology, this is a compact space on which $G$ acts by left  multiplication.
For a reduced word $s$ in $G$, write $\Omega(s)$ for the  cylinder set consisting of all infinite
words in
$\Omega$ that begin with 
$s$. Let $\Omega(1) = \Omega$. Defining a probability measure $\mu$ on  
$\Omega$ amounts to specifying $\mu(\Omega(s))$ in [0,1] for each $s$ in $G$  in such a way
that $\mu(\Omega(1)) = 1$ and
$$\mu(\Omega(s)) = \sum\{\mu(\Omega(sv)) : |v| = 1, |sv| = |s| + 1 \} \ .$$ If
$\mu$ is quasiinvariant under the left action of $G$, and if $p_1, p_2, \ldots, p_n$ are
complex-valued measurable functions on $\Omega$ such that 
$$|p_j(\omega)|^2 = \frac{d\mu(u_j^{-1} \omega)}{d\mu(\omega)}$$ for $\mu-$almost all
$\omega$, then we obtain a unitary representation $\pi$ of $G$ on $L^2(\Omega, \mu)$ by
sending each $u_j$ to the unitary $U_j$ defined on $L^2$ by $$(U_j \xi)(\omega) = p_j(\omega)
\xi(u_j^{-1} \omega) \ .$$ 

The particular type of measure on $\Omega$ that we want to consider here is defined in terms of
positive real numbers $\alpha_+, \alpha_-, \alpha_0, \alpha_1$ satisfying
$$n (\alpha_+ + \alpha_-) = 1 = n \alpha_0 + (n-1) \alpha_1 .$$ Let $\mu(\Omega(u_j)) =
\alpha_+$ and $\mu(\Omega(u_j^{-1})) = \alpha_-$ for each $j$. Once $\mu(\Omega(\cdot))$ has
been defined on words of length $k$, extend the definition to words of length $k+1$ by letting
$\mu(\Omega(sv))$ (where
$|v|= 1 = |sv| - |s|$) be either $\mu(\Omega(s))\alpha_1$ or
$\mu(\Omega(s))\alpha_0$ depending on whether $sv$ ends or doesn't end with a generator
exponent change. (For example, $\mu(\Omega(u_2 u_1 u_2)) = \alpha_+ \alpha_0^2,$ while
$\mu(\Omega(u_2 u_1^{-1} u_2)) = \alpha_+ \alpha_1^2,$ and $\mu(\Omega(u_2 u_1 u_2^{-1})) =
\alpha_+ \alpha_0 \alpha_1.$) It is then straightforward to show that
$$\frac{d\mu(u_j^{-1} \omega)}{d\mu(\omega)} = \left \{ 
\begin{array}{cl}
\alpha_- \alpha_1/\alpha_+& \ \omega \in \Omega(u_i) \ (i \not= j) \\
&\\
\alpha_0& \ \omega \in \Omega(u_i^{-1}) \ (\mbox{any} \ i) \\
&\\
1/\alpha_0& \ \omega \in \Omega(u_j u_i) \ (\mbox{any} \ i) \\
&\\
\alpha_-/(\alpha_+ \alpha_1)& \ \omega \in \Omega(u_j u_i^{-1}) 
\ (i \not= j)
\end{array} \ . \right .$$

Take $\lambda$ in the real interval $(0, \sqrt{n}).$ Define
$\alpha$'s for the construction above by
$$\alpha_+ = \frac{(n-1) \lambda^2}{n (n^2 - \lambda^2)} \ , \ 
\alpha_- = \frac{n-\lambda^2}{n^2 - \lambda^2} \ , \  \alpha_0 =
\frac{\lambda^2}{n^2} \ , \ \alpha_1 = \frac{n - \lambda^2}{n (n - 1)} \ ,$$ and let $\mu$ be
the corresponding measure on  $\Omega.$ Define the functions $p_j$ by
$$p_j(\omega) = \left \{ 
\begin{array}{cl}
(\lambda - n/\lambda)/(n-1)& \ \omega \in
\Omega(u_i) \ (i \not= j) \\
&\\
\lambda/n& \ \omega \in \Omega(u_i^{-1}) \ (\mbox{any} \ i) \\
&\\
n/\lambda& \ \omega \in \Omega(u_j)
\end{array} \ . \right .$$ One checks readily that  $p_j$ is appropriately related to the
Radon-Nikodym derivative of the translate of $\mu$ by $u_j$. Let $\pi$ be the resulting
representation of $G$ on $L^2(\Omega, \mu),$ and let $\xi_0$ be the constant function 1 on
$\Omega.$ Define a state $\phi$ on $G$ by $\phi(s)\ = \ <\pi(s) \xi_0, \xi_0>,$ that is,
$$\phi(s) = \int_\Omega P(s,\omega) \, d\mu(\omega),$$ where $P(\cdot \ , \cdot)$ is the
(nonzero-real-valued) cocycle such that
$P(u_j, \cdot) = p_j$ for each $j$. A routine but tedious calculation, which we omit,
establishes that this state $\phi$ coincides with the state $\phi_{\lambda/n}$ constructed in
Section \ref{eigenstates}. 

The calculation we are leaving out here may be viewed as a mild test of the conjecture set
forth in the introduction to this paper. It is more or less immediate that the sum of the
functions $p_j$ above is the constant function $\lambda,$ in other words that $\xi_0$ is a
$\lambda-$eigenvector for $\pi(u_1 + \ldots + u_n)$, which makes $<\pi(\cdot) \xi_0, \xi_0>$
a $\lambda-$eigenstate for the sum of the generators. Furthermore, Theorem
2.7 in \cite{Spielberg} (see also Theorem 1X in \cite{KuhnSteger}) says that $\pi$ (like all
other boundary representations) is weakly contained in the left regular representation, so one
knows that this state is reduced even before checking that it coincides with $\phi_{\lambda/n}$.
It would appear at first glance that breaking some of the symmetry in the construction above
ought to yield a great many reduced $\lambda-$eigenstates for the sum of the generators. The only
requirements are: (1) a quasiinvariant probability measure $\mu$; and (2) (measurable) complex
functions $z_1, \ldots, z_n$ on $\Omega$ of modulus 1 such that
$$ \sum_{j=1}^n z_j(\omega) \sqrt{d\mu(u_j^{-1}\omega)/d\mu(\omega)} = \lambda$$
for almost every $\omega.$ Condition (2) seems not particularly
onerous; for instance, in the case of two generators, (2) amounts to 
$$|\sqrt{d\mu(u_1^{-1}\omega)/d\mu(\omega)} - \sqrt{d\mu(u_2^{-1}\omega)/d\mu(\omega)}| \leq
\lambda$$ $$ \leq \sqrt{d\mu(u_1^{-1}\omega)/d\mu(\omega)} +
\sqrt{d\mu(u_2^{-1}\omega)/d\mu(\omega)}$$
almost everywhere. This seems to be not so easy to achieve, however. Modest numerical
experimentation instead favors the conjecture that the essential supremum of 
$$|\sqrt{d\mu(u_1^{-1}\omega)/d\mu(\omega)} -
\sqrt{d\mu(u_2^{-1}\omega)/d\mu(\omega)}|$$
is always greater than or equal to the essential infimum of
$$ \sqrt{d\mu(u_1^{-1}\omega)/d\mu(\omega)} +
\sqrt{d\mu(u_2^{-1}\omega)/d\mu(\omega)} \ , $$
with equality only in the situation of the previous paragraph. In any event, questions of
uniqueness of reduced eigenstates plainly have a good deal to do with the behavior of
quasiinvariant measures on $\Omega$.

\

\section{Paucity of eigenstates}\label{paucity}

\

We conjecture that for $|\lambda| \leq \sqrt{n}$, there is only one reduced
$\lambda-$eigenstate for $u_1 + \ldots + u_n - \lambda,$ namely the one exhibited in Section
\ref{eigenstates}. We will prove this below in case $\lambda = \sqrt{n}$ (so by rotating, for
$|\lambda| = \sqrt{n}$).

Henceforth we will work in $C^*_r(G)$, and think of group elements and the group algebra as
acting on $\ell^2(G)$ via the left regular representation. Let 
$$T = \frac{1}{\sqrt{n}} (u_1 + \ldots + u_n) \ ,$$  so $T$ belongs to $L(\ell^2(G)),$ the
algebra of bounded operators on $\ell^2(G).$ Let
$S^+$ be set of reduced words in $G$ beginning with a positive generator power, and let $S^- = G
\setminus S^+$, so $S^-$ consists of 1 together with the words beginning with a negative
generator power. Let $P$ be the orthogonal projection of $\ell^2(G)$ on $\ell^2(S^+)$, and write
$Q = I - P,$ the projection on $\ell^2(S^-).$ 

The state whose uniqueness we are trying to prove is $\phi_{1/\sqrt{n}}$ from Section
\ref{eigenstates}, which we will write simply as $\phi.$ Notice that $\phi(s) = n^{-|s|/2}$ if
$\gamma(s) = 0$ (that is, if $s \in G^+(G^+)^{-1})$ and $\phi(s) = 0$ if $\gamma(s) > 0.$

\

\begin{lem}\label{proto} (a) If $f$ is a state of $L(\ell^2(G))$ such that 
$$f(P) = 1 \ \  \mbox{and} \ \  f((T^* - I) (T - I)) = 0 \ ,$$ 
then $f(s) = \phi(s)$ for all $s$ in $G$.

(b) If instead 
$$f(Q) = 1 \ \  \mbox{and} \ \ f((T - I) (T^* - I)) = 0 \ ,$$ 
then $f(s) = \phi(\sigma(s))$ for all $s$ in $G$, where $\sigma$ is the automorphism 
of $G$ that sends each generator $u_j$ to its inverse.
\end{lem}

\

\noindent{\bf Proof:} (a)  We show first that $f(s) = 0$ for $s$ not in $G^+ (G^+)^{-1}$. Indeed, if
$\gamma(s) > 0,$ then $s$ contains $u_i^{-1} u_j$ for some $i \not= j$. For sufficiently large
$m$, and any $t_1, t_2$ in $G^+$ of length $m$, we can write $t_1^{-1} s t_2$ as 
$v_1^{-1} u_i^{-1} u_j v_2$ with no cancellation, where $v_2, v_1$ end in positive generator
powers. Notice that 
$$v_1^{-1} u_i^{-1} u_j v_2 S^+ \subseteq S^-,$$ 
so $P v_1^{-1} u_i^{-1} u_j v_2 P = 0.$ It follows that
$P(T^*)^m \, s \, T^m P = 0$ for sufficiently large $m$. The hypotheses on $f$ make
$f(P(T^*)^m X T^m P) = f(X)$ for all $X$ in $L(\ell^2(G))$, so $f(s) = 0.$ 

It is now easy to see that $f$ and $\phi$ coincide on $G^+ (G^+)^{-1}$. Namely for any $s$ in
$G^+ (G^+)^{-1}$, including $s = 1,$ we have
$$f(s u_j^{-1}) = f(s u_j^{-1} T) = \frac{1}{\sqrt{n}} f(s)$$
because $\gamma(su_j^{-1} u_i) > 0$ for $i \not= j.$ Likewise
$$f(u_j s) = \frac{1}{\sqrt{n}}f(s).$$ Since $f(1) = 1,$ iteration gives $f(s) = n^{-|s|/2}$ when
$\gamma(s) = 0.$

(b) The slight asymmetry between (a) and (b) that results from our putting $1$ in $S^-$ 
rather than in $S^+$ is harmless because 
$$v_1^{-1} u_i^{-1} u_j v_2 \in S^- \setminus \{1\}$$
in the argument above for (a).

\

The next lemma records two simple observations.

\

\begin{lem}\label{obs} (a) $QTT^*Q = Q$; 
\ \ \ \ \ (b) $P u_i^{-1} u_j T^* Q = 0$ for $i \not= j$.
\end{lem}

\

\noindent{\bf Proof:} (a) Notice that
$$T T^*  = \frac{1}{n} \left( n I + 
\sum_{i \not= j} u_i u_j^{-1} \right) .$$
For $i \not= j,$ we have $Q u_i u_j^{-1} Q = 0.$

(b) Plainly, $P u_i^{-1} Q = 0$ and furthermore $P u_i^{-1} u_j u_k^{-1} Q = 0$
provided $k \not= j.$ Thus
$$P u_i^{-1} u_j T^* Q = 
\frac{1}{\sqrt{n}} \sum_{k=1}^n P u_i^{-1} u_j u_k^{-1} Q = 0.$$

\

\begin{thm}\label{unique} The only reduced $\sqrt{n}-$eigenstate of $u_1 + \ldots u_n$ is $\phi,$
given by 
$$\phi(s) = \left\{\begin{array}{cl} 
n^{-|s|/2} & s \in G^+ (G^+)^{-1} \\
0 & \mbox{else} \end{array} \right . \ . $$
\end{thm}

\

\noindent{\bf Proof:} Let $\psi$ be such an eigenstate. Identifying $\psi$ with a state of $C^*_r(G)$ and
then extending to a state $f$ on all of $L(\ell^2(G))$, we have $f((T^* - I) (T - I)) = 0.$ In
light of Lemma \ref{proto}(a), we need only show $f(P) = 1$ (that is, $f(Q) = 0$) in order to
prove the theorem. Suppose that $f(Q) > 0.$ Let $g$ be the state of $L(\ell^2(G))$ defined by
$g(X) = f(QXQ)/f(Q).$ Then $g(Q) = 1,$ and furthermore
$$g((I - T) (I - T^*)) = f(Q)^{-1} (f(Q) - f(QTQ) - f(QT^*Q) + f(QTT^*Q)).$$ Because $QT^*Q =
T^*Q$, and because $I - T^*$ is in the right kernel of $f$, we have
$f(QT^*Q) = f(Q)$. Taking adjoints gives $f(QTQ) = f(Q).$ Lemma \ref{obs}(a) now makes
$g((I - T) (I - T^*)) = 0.$ It follows from Lemma \ref{proto}(b) that 
$g(u_i^{-1} u_j) = 1/n$ for $i \not= j$. Consider now $f(u_i^{-1} u_j).$ We have
$$f(P u_i^{-1} u_j Q) = f(P u_i^{-1} u_j T^*Q) = 0$$ by Lemma \ref{obs}(b). (We have already
shown that $f((Q - QT) (Q - T^*Q)) = 0,$ so $f(X Q) = f(X T^* Q)$ for all $X$ in
$L(\ell^2(G))$.) Taking adjoints and swapping $i$ and $j$ shows that $f(Q u_i^{-1} u_j P) =
0.$ Plainly $P u_i^{-1} u_j P = 0$, so we have 
$$f(u_i^{-1} u_j) = f(Q u_i^{-1} u_j Q) = f(Q)
g(u_i^{-1} u_j)  =
\frac{f(Q)}{n}$$ for $i \not= j.$ Now $f(T^*T) = f(I) = 1$ because $I - T$ is in the left
kernel of $f$. This forces the sum over unequal $i$ and $j$ of $f(u_i^{-1} u_j)$ to vanish,
but by what we have just shown, this sum is $(n-1) f(Q).$ Thus $f(Q)$ must after all be 0.

\ 

We conclude this paper with a look at the states of $C^*_r(G)$ that have a polynomial in one
of the generators in their left kernel.  The simplest case is $u_1 - z,$ where $z$ is a complex
scalar of modulus 1. Let
$G_1$ be the subgroup of $G$ generated by
$u_1,$  and let $\chi_1$ be the characteristic function of $G_1$. It is easily
checked that 
$$\chi_1(s) = \lim_{k \rightarrow \infty} \frac{1}{k} <s \xi_k, \xi_k>$$
for every $s$ in $G$, where  $\xi_k$ is the characteristic function of $\{u_1, u_1^2, \ldots,
u_1^k\}$, so $\chi_1$ is a reduced state of $G$. It is plainly also a $1-$eigenstate for $u_1.$
To obtain $z-$eigenstates (for $|z| = 1$), precede $\chi_1$ by the automorphism of the group
algebra that sends $u_1$ to $zu_1$ and fixes the other generators. This gives the reduced
$z-$eigenstate $\chi_z,$ which satisfies $\chi_z(u_1^j) = z^j$ and $\chi_z(s) = 0$ for $s$ in
$G \setminus G_1$.  Let us write $\pi_z$ for the unitary representation of $G$ constructed
from $\chi_z$. 

The representations $\pi_z$ were first considered by H. Yoshizawa
in \cite{Yoshizawa}, where it is shown that these representations are irreducible,
unitarily inequivalent to one another, and weakly contained in the left regular representation.
These facts also follow from the lemma below, of course. Notice, by the way, that
$\pi_1$ is the representation of $G$ that comes from its left action on $G/G_1.$

\

\begin{lem}\label{chi} For $|z| = 1,$ the only reduced state $\psi$ of $G$ such that $\psi(u_1) =
z$ is $\chi_z.$
\end{lem}

\

\noindent{\bf Proof:} We may take $z = 1$ without loss of generality. The Cauchy-Schwarz
inequality shows that
$\psi(u_1 s) = \psi(s u_1) = \psi(s)$ for every $s$ in $G$, so $\psi$ is identically 1 on $G_1$,
and constant on each double coset $G_1 s G_1$. If $s \notin G_1,$ then 0 belongs to the
norm-closed convex hull of 
$G_1 s G_1$ in $C^*_r(G)$. (This follows from Theorem IV J in \cite{AkemannOstrand} because
the double coset contains an infinite free subset of $G$, and also from the averaging result
in \cite{Archbold}.) Thus $\psi(s) = 0.$

\

The following theorem (whose proof was kindly supplied by the referee of this paper in place of
a more involved argument) shows that only finitely many reduced pure states of
$G$ can have a polynomial in one of the generators in their left kernel, namely (if the generator is
$u_1$) the $\chi_z$'s for unimodular zeros $z$ of the polynomial.

\

\begin{thm}\label{poly} Let $p$ be a polynomial with complex coefficients, and let $\psi$ be a
reduced state of $G$ such that $\psi(|p(u_1)|^2) = 0.$ Then $\psi$ must be a convex combination
of $\chi_z$'s for $z$'s among the zeros of $p$ on the unit circle. Said another way, the only
representations $\pi$ of $G$ weakly contained in the left regular representation such that
$\pi(p(u_1))$ has a nonzero kernel are direct sums of $\pi_z$'s (with $p(z) = 0$).
\end{thm}

\

\noindent{\bf Proof:} Let $H, \pi, \xi$ be the Hilbert space, unitary representation, and
unit cyclic vector constructed from $\pi,$ so $\psi(s) \ =  \ <\pi(s) \xi, \xi>.$ Thus,
$\pi(p(u_1)) \xi = 0.$ Restricting to the $C^*-$subalgebra $A$ of $C^*_r(G)$generated by $u_1$
--- which we identify in the usual way with the algebra of continuous complex functions on
the unit circle --- we find that $\psi|_A$ is a convex combination of point masses at zeros of
$p$ on the circle. By considering $\pi(g(u_1)) \xi$ for continuous functions $g$ on the circle
that take the value 1 at a particular zero of $p$ and vanish at all the other zeros that $p$ has
on the circle, we obtain nonzero vectors $\xi_1, \ldots, \xi_k$ in $H$ and
distinct modulus-one zeros $z_1, \ldots, z_k$ of $p$ such that
$$\pi(f(u_1)) \xi \ = \ \sum_{j=1}^k f(z_j) \xi_j$$
for all continuous functions $f$ on the circle. By construction (or for other reasons), the
vectors $\xi_j$ are pairwise orthogonal and sum to $\xi.$ (Only the zeros of $p$ that
contribute positively to $\psi|_A$ as a convex combination of point masses are counted among the
$z_j$'s; we ignore the other zeros of $p$.) For $j = 1, \ldots, k,$ let $\alpha_j = ||
\xi_j||^2$ and let $H_j$ be the closed subspace of $H$ spanned by $\pi(G) \xi_j.$ 
Since $<\pi(u_1) \xi_j, \xi_j> = <\pi(u_1) \xi, \xi_j> = \alpha_j z_j,$ the reduced state
$\alpha_j^{-1} <\pi(\cdot) \xi_j, \xi_j>$ must by \ref{chi} be $\chi_{z_j}$ and the
representation obtained by restricting $\pi$ to $H_j$ must be the Yoshizawa
representation $\pi_{z_j}$. These subrepresentations are therefore irreducible and unitarily
inequivalent. It follows that the subspaces $H_j$ are orthogonal to one another. (If
$P_i$ is the orthogonal projection of $H$ on $H_i$, then the restriction of $P_i$
to $H_j$ interwines the restrictions of $\pi$ to $H_j$ and $H_i$.) The
orthogonality of the $H_j$'s makes
$$<\pi(s) \xi, \xi> \ = \ \sum_{j=1}^k <\pi(s) \xi_j, \xi_j> \ = \ \sum_{j=1}^k \alpha_j
\chi_{z_j}(s)$$ for all $s$ in $G$.

\

\

\

\


\begin{thebibliography}{}
\bibitem[AO]{AkemannOstrand}  C. A. Akemann and P. A. Ostrand, {\em Computing norms in group
C*-algebras}, Amer. J. Math. {\textbf 98} (1976), 1015 -- 1047. 
\bibitem[A]{Archbold}  R. J. Archbold, {\em A mean ergodic theorem associated with the free
group on two generators}, J. London Math. Soc. {\textbf 13} (1976), 339 -- 345.
\bibitem[DMFT]{DeMicheleFigaTalamanca}  L. De Michele and A. Fig\`{a}-Talamanca, {\em Positive
definite functions on free groups}, Amer. J. Math. {\textbf 102} (1980), 503 -- 509.
\bibitem[FTP]{FigaTalamancaPicardello} A. Fig\`{a}-Talamanca and M. A. Picardello, {\em
Harmonic Analysis on Free Groups}, Lecture Notes in Pure and Applied Mathematics, {\textbf 87},
Marcel Dekker, New York, 1983.
\bibitem[FTS]{FigaTalamancaSteger} A. Fig\`{a}-Talamanca and T. Steger, {\em Harmonic analysis for
anisotropic random walks on homogeneous trees}, Mem. Amer. Math. Soc. {\textbf 110} (1994), \# 531.
\bibitem[H]{Haagerup} U. Haagerup, {\em An example of a non nuclear C*-algebra which has the
metric approximation property}, Invent. Math. {\textbf 50} (1979), 279 -- 293.
\bibitem[KS]{KuhnSteger} G. Kuhn and T. Steger, {\em More irreducible boundary
representations}, Duke Math. J. {\textbf 82} (1996), 381 -- 436.
\bibitem[L]{Linnell} P. A. Linnell, {\em Division rings and group von Neumann algebras}, Forum
Math. {\textbf 5} (1993), 561 -- 576.
\bibitem[S]{Spielberg} J. Spielberg, {\em Free-product groups, Cuntz-Krieger algebras, and
covariant maps}, International J. Math. {\textbf 2} (1991), 457 -- 476.
\bibitem[Y]{Yoshizawa} H. Yoshizawa, {\em Some remarks on unitary representations of the free
group}, Osaka J. Math. {\textbf 3} (1951), 55 --- 63.
\end{thebibliography}
\end{document}